# SCHUR FUNCTIONS AND ORTHOGONAL POLYNOMIALS ON THE UNIT CIRCLE

Ferenc Pintér and Paul Nevai

The Ohio State University

ABSTRACT. We apply a theorem of Geronimus to derive some new formulas connecting Schur functions with orthogonal polynomials on the unit circle. The applications include the description of the associated measures and a short proof of Boyd's result about Schur functions. We also give a simple proof for the above mentioned theorem of Geronimus.

## 1. Schur functions

In what follows we adopt the following notations: $\mathbb{D} \stackrel{\text{def}}{=} \{z \in \mathbb{C} : |z| < 1\}$, $\overline{\mathbb{D}}$ stands for the closure of $\mathbb{D}$, and $\mathbb{Z}^+ \stackrel{\text{def}}{=} \mathbb{N} \cup \{0\}$. In addition, $\mathcal{B}$ denotes the set of Schur functions, namely,

$$\mathcal{B} \stackrel{\text{def}}{=} \{f \,:\, f \text{ is analytic and } |f| < 1 \text{ in } \mathbb{D}\}.$$

Similarly, $\mathcal{C}$ stands for the set of Carathéodory functions, that is,

$$\mathcal{C} \stackrel{\text{def}}{=} \{F \,:\, F \text{ is analytic and } \Re F > 0 \text{ in } \mathbb{D}, F(0) = 1\}.$$

If $f$ is a Schur function, then

$$F(z) = \frac{1 + z\,f(z)}{1 - z\,f(z)} \qquad (1)$$

is a Carathéodory function, and, vice versa, if $F$ is a Carathéodory function, then

$$f(z) = \frac{1}{z}\,\frac{F(z) - 1}{F(z) + 1} \qquad (2)$$

1991 *Mathematics Subject Classification*. 42C05, 30B70.

This material is based upon work supported by the National Science Foundation under Grant No. DMS–940577.

Typeset by $\mathcal{A}_{\mathcal{M}}\mathcal{S}$-TeX

is a Schur function. The starting point of Schur's algorithm is a function $f \in \mathcal{B}$ which is used to generate a sequence of functions in $\mathcal{B}$ in the following way

$$f_0 \stackrel{\text{def}}{=} f \quad \& \quad \gamma_n \stackrel{\text{def}}{=} f_n(0) \quad \& \quad f_{n+1}(z) \stackrel{\text{def}}{=} \frac{1}{z} \frac{f_n(z) - \gamma_n}{1 - \overline{\gamma}_n f_n(z)}, \quad n \in \mathbb{Z}^+. \quad (3)$$

The $\gamma_n$'s are called the Schur parameters corresponding to $f \in \mathcal{B}$. Iteratively applying this definition, two sequences of polynomials $\{A_n\}_0^\infty$ and $\{B_n\}_0^\infty$ are formed to yield

$$f_{n+1}(z) = \frac{1}{z} \frac{f(z) B_n(z) - A_n(z)}{B_n^*(z) - f(z) A_n^*(z)}, \quad (4)$$

or, equivalently,

$$f(z) = \frac{A_n(z) + z B_n^*(z) f_{n+1}(z)}{B_n(z) + z A_n^*(z) f_{n+1}(z)}, \quad (5)$$

where the reversed $*$-polynomial of a polynomial $\varrho_n$ of degree $n$ is defined by $\varrho_n^*(z) \stackrel{\text{def}}{=} z^n \overline{\varrho}_n(z^{-1})$. The polynomial sequences $\{A_n\}_0^\infty$ and $\{B_n\}_0^\infty$ can be generated by the following matrix recurrences (cf. [1, § 1, formulas (1)–(6), p. 145])

$$A_0 \stackrel{\text{def}}{=} \gamma_0 \quad \& \quad B_0 \stackrel{\text{def}}{=} 1 \quad \&$$

$$\begin{pmatrix} A_n^* & B_n^* \\ -B_n & -A_n \end{pmatrix} = \begin{pmatrix} z & -\overline{\gamma}_n \\ -z \gamma_n & 1 \end{pmatrix} \begin{pmatrix} A_{n-1}^* & B_{n-1}^* \\ -B_{n-1} & -A_{n-1} \end{pmatrix}, \quad n \in \mathbb{N}. \quad (6)$$

If

$$\omega_n \stackrel{\text{def}}{=} \prod_{k=0}^n \left(1 - |\gamma_k|^2\right), \quad (7)$$

then taking the determinants in (6), we obtain

$$B_n(z) B_n^*(z) - A_n(z) A_n^*(z) = z^n \omega_n, \quad n \in \mathbb{Z}^+. \quad (8)$$

In connection to Schur's classical result, the following theorem is known (cf. [12, the theorem in § 3.I, p. 211, and § 14, p. 137–141]).

**Theorem 1.** *If $f \in \mathcal{B}$ then $\lim_{n \to \infty} A_n/B_n = f$ locally uniformly in $\mathbb{D}$, where $\{A_n\}_0^\infty$ and $\{B_n\}_0^\infty$ are generated by* (6).

*Proof.* Our proof is based on [1, § 1, p. 145–146]. By (5) and (8) we can write

$$f(z) - \frac{A_n(z)}{B_n(z)} = \frac{A_n(z) + z B_n^*(z) f_{n+1}(z)}{B_n(z) + z A_n^*(z) f_{n+1}(z)} - \frac{A_n(z)}{B_n(z)}$$

$$= \frac{z^{n+1} f_{n+1}(z) \omega_n}{B_n(z) (B_n(z) + z A_n^*(z) f_{n+1}(z))}. \quad (9)$$

Applying (5) with $n-1$ instead of $n$ and replacing $f_n(z)$ in the numerator by $(\gamma_n + z\, f_{n+1}(z))/(1 + z\,\overline{\gamma}_n\, f_{n+1}(z))$ (cf. (3)), we obtain

$$f(z) = \frac{A_{n-1}(z) + z\,\gamma_n\, B^*_{n-1}(z) + z\left(\overline{\gamma}_n\, A_{n-1}(z) + z\, B^*_{n-1}(z)\right) f_{n+1}(z)}{\left(B_{n-1}(z) + z\, A^*_{n-1}(z)\, f_n(z)\right)\left(1 + z\,\overline{\gamma}_n\, f_{n+1}(z)\right)}. \qquad (10)$$

By (6), the numerators in (10) and (5) are identical so that the denominators are identical as well. In other words,

$$B_n(z) + z\, A^*_n(z)\, f_{n+1}(z) = \left(B_{n-1}(z) + z\, A^*_{n-1}(z)\, f_n(z)\right)\left(1 + z\,\overline{\gamma}_n\, f_{n+1}(z)\right), \qquad (11)$$

from which

$$B_n(z) + z\, A^*_n(z)\, f_{n+1}(z) = \prod_{k=0}^{n}(1 + z\,\overline{\gamma}_k\, f_{k+1}(z)) \neq 0, \qquad z \in \mathbb{D}, \qquad (12)$$

follows immediately. By (8), $|B_n(z)| > |A_n(z)|$ for $|z| = 1$ so that, by Rouché's theorem, $B_n(z) + z\, A^*_n(z)\, f_{n+1}(z)$ has the same number of roots in $|z| < 1$ as $B_n(z)$. Thus, by (12), $B_n(z)$ has no roots in $|z| \leqslant 1$, and, by the maximum principle, $|A_n(z)/B_n(z)| < 1$ on $|z| \leqslant 1$. Considering $n$ fixed and letting $z \to 0$, we conclude from (9) that the Taylor series expansion of $f$ and that of $A_n/B_n$ coincide up to the $n$th order, namely, $f(z) = \sum_{k=0}^{\infty} u_k z^k$ and $(A_n/B_n)(z) = \sum_{k=0}^{\infty} s_k^{(n)} z^k$, where $s_k^{(n)} = u_k$ for $0 \leqslant k \leqslant n$. Let $0 < \epsilon < 1$ and $|z| \leqslant 1 - \epsilon$. Then, by Cauchy's integral formula for $H^\infty$,

$$\left| f(z) - \frac{A_n(z)}{B_n(z)} \right| \leqslant \left| \sum_{k=n+1}^{\infty} u_k z^k \right| + \left| \sum_{k=n+1}^{\infty} s_k^{(n)} z^k \right|$$

$$\leqslant \sum_{k=n+1}^{\infty} \left| \frac{1}{2\pi} \oint_{|\zeta|=1} \frac{f(\zeta)}{\zeta^{k+1}}\, d\zeta \right| |z^k| + \sum_{k=n+1}^{\infty} \left| \frac{1}{2\pi} \oint_{|\zeta|=1} \frac{(A_n/B_n)(\zeta)}{\zeta^{k+1}}\, d\zeta \right| |z^k|.$$

This can be estimated from above by $2\,(1-\epsilon)^{n+1}/\epsilon$, implying uniform convergence in $|z| \leqslant 1 - \epsilon$ as $n \to \infty$. $\square$

The following concepts will play key roles as we make connection from Schur functions to Szegő's theory.

**Definition 2.** We call $f \in \mathcal{B}$ an extreme point of $\mathcal{B}$ if $f$ is not a proper convex combination of two distinct elements of $\mathcal{B}$.

**Definition 3.** If $g$ is measurable and $\log^+ |g(e^{i\theta})|$ is integrable in $[0, 2\pi]$ (as a function of $\theta$), then we define the outer function $\mathcal{G}(g)$ of $g$ by

$$\mathcal{G}(g; z) \stackrel{\mathrm{def}}{=} \exp\left\{ \frac{1}{2\pi} \int_0^{2\pi} \frac{e^{i\theta} + z}{e^{i\theta} - z} \log |g(e^{i\theta})|\, d\theta \right\}, \qquad |z| < 1. \qquad (13)$$

The following theorem characterizes extreme points in several ways.

**Theorem 4.** *Let $f \in \mathcal{B}$. Then the following statements are equivalent.*
(1) *$f$ is an extreme point of $\mathcal{B}$,*
(2) $\mathcal{G}\left(1 - |f|^2; 0\right) = 0$,
(3) $\int_0^{2\pi} \log\left[1 - |f(e^{i\theta})|\right] d\theta = -\infty$,
(4) $\prod_{k=0}^{\infty} \left(1 - |\gamma_k|^2\right) = 0$,
(5) $\sum_{k=0}^{\infty} |\gamma_k|^2 = \infty$.

*Proof.*

(1)$\Longleftrightarrow$(2): cf. [8, Theorem, p. 138–139].

(2)$\Longleftrightarrow$(4): cf. [1, § 3, Lemma, p. 146–147], and, in a different way, [6, § 2, Theorem, p. 460–461].

(2)$\Longleftrightarrow$(3): since $\mathcal{G}\left(1 - |f|^2\right) = \mathcal{G}(1 - |f|)\,\mathcal{G}(1 + |f|)$, and $1 \leqslant \mathcal{G}(1 + |f|) \leqslant 2$.

(4)$\Longleftrightarrow$(5): true in general. $\square$

## 2. Orthogonal polynomials on the unit circle

Let $\mathbb{T} \stackrel{\text{def}}{=} \{z \in \mathbb{C} : |z| = 1\}$. If $\sigma$ is a positive measure on $[0, 2\pi)$ and $f$ is a function on $\mathbb{T}$, we define $\int_{\mathbb{T}} f\, d\sigma \stackrel{\text{def}}{=} \int_0^{2\pi} f(e^{i\theta})\, d\sigma(\theta)$ when the right-hand side exists. Orthogonal polynomials $\{\varphi_n\}_0^{\infty}$ on the unit circle are defined by

$$\frac{1}{2\pi} \int_{\mathbb{T}} \varphi_n(\sigma, z) \overline{\varphi_m(\sigma, z)}\, d\sigma = \delta_{m,n}, \qquad m, n \in \mathbb{Z}^+, \tag{14}$$

where $\sigma$ is a positive measure on $[0, 2\pi)$ with infinite support. If $\varphi_n(z) = \kappa_n z^n +$ lower degree terms, $\kappa_n > 0$, then $\Phi_n(z) \stackrel{\text{def}}{=} \kappa_n^{-1} \varphi_n(z)$ are called monic orthogonal polynomials. They satisfy the Szegő recurrences

$$\Phi_0 \stackrel{\text{def}}{=} 1 \quad \& \quad \Psi_0 \stackrel{\text{def}}{=} 1 \quad \&$$

$$\begin{pmatrix} \Phi_{n+1} & \Psi_{n+1} \\ \Phi_{n+1}^* & -\Psi_{n+1}^* \end{pmatrix} = \begin{pmatrix} z & -\overline{a}_n \\ -z\, a_n & 1 \end{pmatrix} \begin{pmatrix} \Phi_n & \Psi_n \\ \Phi_n^* & -\Psi_n^* \end{pmatrix}, \qquad n \in \mathbb{Z}^+, \tag{15}$$

where $a_n = -\overline{\Phi_{n+1}(0)}$, (cf. [4, formulas (1.2) and (1.2'), p. 6]). The monic second kind orthogonal polynomials $\{\Psi_n\}_0^{\infty}$ are determined by replacing $a_n$ with $-a_n$ in the recurrences for $\{\Phi_n\}_0^{\infty}$ and $\{\Phi_n^*\}_0^{\infty}$. We can connect the leading coefficients $\{\kappa_n\}_0^{\infty}$ to the reflection coefficients $\{a_n\}_0^{\infty}$ via

$$\sum_{k=0}^{n} |\varphi_k(0)|^2 = \kappa_n^2, \qquad n \in \mathbb{Z}^+, \tag{16}$$

and

$$\frac{\kappa_n^2}{\kappa_{n+1}^2} = 1 - |a_n|^2, \qquad n \in \mathbb{Z}^+, \tag{17}$$

(cf. [4, formula (1.5), p. 7, or formula (1.9), p. 9]). It is also a well known fact, which follows from (15) by induction, that

$$\Psi_n(z) = \frac{1}{2\pi} \int_0^{2\pi} \frac{e^{i\theta} + z}{e^{i\theta} - z} \left(\Phi_n(e^{i\theta}) - \Phi_n(z)\right) d\sigma(\theta), \qquad n \in \mathbb{N}. \tag{18}$$

Since $\Phi_n^*$ has all its zeros outside the unit circle (cf. [4, formula (1.12), p. 9]), (18) implies that

$$\lim_{n \to \infty} \frac{\Psi_n^*(z)}{\Phi_n^*(z)} = \frac{1}{2\pi} \int_0^{2\pi} \frac{e^{i\theta} + z}{e^{i\theta} - z} d\sigma(\theta), \tag{19}$$

locally uniformly in $\mathbb{D}$ (cf. [4, formula (1.16'), p. 11]). In the theory of orthogonal polynomials, Szegő's theory plays a distinct role. Its scope covers the case when one of the following equivalent conditions (cf. the Szegő condition) holds.

$$\log \sigma' \in \mathrm{L}^1(\mathbb{T}) \quad \Longleftrightarrow \quad \sum_{k=0}^{\infty} |\varphi_k(\sigma, 0)|^2 < \infty \quad \Longleftrightarrow \quad \lim_{n \to \infty} \kappa_n(\sigma) < \infty \quad \Longleftrightarrow$$

$$\overline{\mathbb{P}} \neq \mathrm{L}^2(\sigma, \mathbb{T}) \quad \Longleftrightarrow \quad \prod_{k=0}^{\infty} \left(1 - |a_k|^2\right) > 0 \quad \Longleftrightarrow \quad \sum_{k=0}^{\infty} |a_k|^2 < \infty,$$

where $\mathbb{P}$ denotes the set of complex algebraic polynomials (cf. [7, Theorem 3.1(a), p. 44, and Theorem 3.3.(a), p. 48] and [4, Section 1.1.14, p. 14–18]).

## 3. Relations between Schur functions and orthogonal polynomials on the unit circle

The basic connection between Schur functions and orthogonal polynomials on the unit circle is established by Herglotz's theorem which says that every Carathéodory function $F$ with $F(0) = 1$ has the representation

$$F(z) = \frac{1}{2\pi} \int_0^{2\pi} \frac{e^{i\theta} + z}{e^{i\theta} - z} d\sigma(\theta),$$

where $\sigma$ is a unique positive Borel measure on the unit circle (cf. [8, exercise 10, p. 40, and Theorem, p. 66–67]). The reverse statement is obviously true. To recover the measure explicitly from $F$, one can use the following Stieltjes inversion formula

$$\frac{\sigma(\theta+) + \sigma(\theta-)}{2} = \mathrm{const} + \lim_{r \to 1-} \int_0^{\theta} \Re F(re^{it}) \, dt.$$

In what follows we assume this $\sigma$ to be the measure of orthogonality. Then we can simultaneously use the concepts introduced in sections 1 and 2, related through the above equations and (2). Geronimus proved that $\gamma_n = a_n$ for $n \in \mathbb{Z}^+$ (cf. [3, Theorem IX, 2°, formula (44), p. 110–112]). However, a similarly important connection exists between the polynomial sequences $\{A_n\}_0^{\infty}$, $\{B_n\}_0^{\infty}$ and $\{\Phi_n\}_0^{\infty}$, $\{\Psi_n\}_0^{\infty}$. This relationship is described by the following theorem.

**Theorem 5.** *With the notations of the previous sections we have for* $n \in \mathbb{Z}^+$

$$A_n(z) = \frac{\Psi^*_{n+1}(z) - \Phi^*_{n+1}(z)}{2z}, \qquad B_n(z) = \frac{\Psi^*_{n+1}(z) + \Phi^*_{n+1}(z)}{2}, \qquad (20)$$

$$A^*_n(z) = \frac{\Psi_{n+1}(z) - \Phi_{n+1}(z)}{2}, \qquad B^*_n(z) = \frac{\Psi_{n+1}(z) + \Phi_{n+1}(z)}{2z}, \qquad (21)$$

*or, in a matrix product form,*

$$\begin{pmatrix} A^*_n & B^*_n \\ -B_n & -A_n \end{pmatrix} = \begin{pmatrix} \Phi_{n+1} & \Psi_{n+1} \\ \Phi^*_{n+1} & -\Psi^*_{n+1} \end{pmatrix} \begin{pmatrix} -1/2 & 1/(2z) \\ 1/2 & 1/(2z) \end{pmatrix}. \qquad (22)$$

*Proof.* The right-hand sides in (20) and (21) are polynomials of order $n$ by means of the following forms

$$\Phi_n(z) = z^n + \cdots + (-\overline{a}_{n-1}),$$
$$\Phi^*_n(z) = (-a_{n-1})z^n + \cdots + 1,$$
$$\Psi_n(z) = z^n + \cdots + \overline{a}_{n-1},$$
$$\Psi^*_n(z) = a_{n-1}z^n + \cdots + 1.$$

Now the proof is an immediate consequence of Geronimus' above mentioned theorem, using (6) and (15) with induction. □

*Remark 6.* Embedded in his proof of Geronimus' theorem, Khrushchev proved half of Theorem 5, namely, that $\Phi^*_n(z) = B_{n-1}(z) - z A_{n-1}(z)$ (cf. [9, formulas (2.6) and (2.10), p. 188–189]). One can also derive (20) and (21), following Golinskii's proof of Geronimus' theorem, see [6, § 2, p. 459–460]. Indeed, the concluding line of his proof can be written as

$$f_n(z) = \frac{1}{z} \frac{z f(z)(\Phi^*_n(z) + \Psi^*_n(z)) + \Phi^*_n(z) - \Psi^*_n(z)}{z f(z)(\Phi_n(z) - \Psi_n(z)) + \Phi_n(z) + \Psi_n(z)}. \qquad (23)$$

Now (20) and (21) are suggested by a comparison of (23) with (4). To complete the proof of Theorem 5, we still need induction, since the uniqueness of the sequences $\{A_n\}_0^\infty$ and $\{B_n\}_0^\infty$ is not guaranteed.

*Remark 7.* Formula (8) along with (22) and (17) imply that

$$\varphi_n(z)\psi^*_n(z) + \psi_n(z)\varphi^*_n(z) = 2z^n, \qquad n \in \mathbb{Z}^+, \qquad (24)$$

(cf. [4, formula (1.17), p. 11–12]).

*Remark 8, A short proof of Geronimus' theorem.* The proof of Theorem 5 also motivates Geronimus' theorem by underlining the similarity of the matrices involved in (6) and (15). There is yet another way to motivate the equations in (20). Compare the convergence results on $\lim_{n\to\infty} \Psi^*_n/\Phi^*_n$ and $\lim_{n\to\infty} A_n/B_n$. The right-hand side of (19) is a Carathéodory function ($\stackrel{\text{def}}{=} F$), but, via (2) the convergence can be transferred to a Schur function ($\stackrel{\text{def}}{=} f$) as well. Considering this $f$ as the initial Schur function to generate the sequence $\{f_n\}_0^\infty \subset \mathcal{B}$, we have

$$\lim_{n\to\infty} \frac{\Psi_n^*(z) - \Phi_n^*(z)}{z\left(\Psi_n^*(z) + \Phi_n^*(z)\right)} = f(z) = \lim_{n\to\infty} \frac{A_n(z)}{B_n(z)}, \qquad (25)$$

which suggests (20) up to some questions about the right degrees of the polynomials above. Assuming (20) and (21) now, we can use (4) to guess that

$$f_n(z) = \frac{F(z)\,\Phi_n^*(z) - \Psi_n^*(z)}{z\left(F(z)\,\Phi_n(z) + \Psi_n(z)\right)}. \qquad (26)$$

The latter, however, can be proven by induction, using (3) and (15), depending on the following relation

$$\lim_{z\to 0} \frac{F(z)\,\Phi_n^*(z) - \Psi_n^*(z)}{z\left(F(z)\,\Phi_n(z) + \Psi_n(z)\right)} = a_n, \qquad (27)$$

which is true by virtue of [4, formulas (1.3) and (1.4), p. 6–7, and formulas (1.15') and (1.16'), p. 10–11]. Let us also notice that by proving (26) and (27), we have also proved Geronimus' theorem. $\square$

*Remark 9.* Favard's theorem on the unit circle states that for an arbitrary sequence $\{a_n\}_0^\infty$, $|a_n| < 1$, there exists a measure $\sigma$ (unique if $\sigma/2\pi$ is a probability measure), such that $a_n = -\overline{\Phi_{n+1}(\sigma,0)}$, $n \in \mathbb{Z}^+$, where $\{\Phi_n(\sigma)\}_0^\infty$ are the polynomials orthogonal with respect to $\sigma$. An analogue of this question regarding Schur functions is, if for an arbitrary sequence $\{\gamma_n\}_0^\infty$, $|\gamma_n| < 1$, there exists $f \in \mathcal{B}$ such that $f_n(0) = \gamma_n$ for $n \in \mathbb{Z}^+$. By Geronimus' theorem these two statements are equivalent. For the uniqueness part of the latter one we need Geronimus' theorem along with (19) only. In Wall's formulation of Schur's algorithm, the fractions $A_n/B_n$ and $zB_n^*/zA_n^*$ are the even and odd order convergents, respectively, of the continued fraction

$$\gamma_0 + \frac{\left(1 - |\gamma_0|^2\right)z}{\overline{\gamma}_0\,z} + \frac{1}{\gamma_1} + \frac{\left(1 - |\gamma_1|^2\right)z}{\overline{\gamma}_1\,z} + \ldots, \qquad (28)$$

(cf. [13, formula (2.1), p. 110], [14, Theorem 77.1, p. 285], and [1, p. 146]),[1] which may be used to give an alternative proof of the existence part of the above quoted theorems.

Boyd's result concerning Schur functions (cf. [1, 4. Theorem, p. 147]) is a consequence of our Theorem 5.

**Theorem 10.** *Suppose that $f \in \mathcal{B}$ is not an extreme point of $\mathcal{B}$. Then $A_n \to a$ and $B_n \to b$ locally uniformly in $\mathbb{D}$ as $n \to \infty$, with $a = b f$. Also $A_n^* \to 0$ and $B_n^* \to 0$ locally uniformly in $\mathbb{D}$.*

*Proof.* To guarantee the existence of $a$ and $b$, and the convergence of $A_n^*$ and $B_n^*$, we apply Theorem 5 first to transform the above statements to statements about orthogonal polynomials on the unit circle. Then we apply Szegő's convergence results [7, formulas (3), (4) and (5) in Theorem 3.4(a), p. 50–51]. To assure that the conditions of the quoted theorem hold, we use $a_n = \gamma_n$ for $n \in \mathbb{Z}^+$ again, and refer to our characterizations of the extreme points of $\mathcal{B}$ and the Szegő condition, respectively. To verify $a = b\,f$, we use Theorem 1. □

---

[1] We thank David W. Boyd for pointing us out that H. S. Wall used a somewhat different definition of Schur functions (cf. [14, p. 285–288, in particular, the definition of $f_1$ on p. 287]) which explains why (28) looks a little different from H. S. Wall's formulas.

## 4. The associated measures and polynomials

To develop the concept of the associated polynomials on the real line, let us recall the continued fractions representation of orthogonal polynomials (cf. [2, p. 85–87]). We call the following continued fraction a Jacobi type or simply J-fraction if the $c_n$'s are real, and the $\lambda_n$'s are positive in

$$\frac{\lambda_1}{x - c_1} - \frac{\lambda_2}{x - c_2} - \frac{\lambda_3}{x - c_3} - \ldots \quad . \tag{29}$$

The $n$th partial denominators $P_n(x)$ of (29) satisfy

$$P_{-1} \stackrel{\text{def}}{=} 0 \quad \& \quad P_0 \stackrel{\text{def}}{=} 1 \quad \&$$

$$P_n(x) = (x - c_n)\,P_{n-1}(x) - \lambda_n\,P_{n-2}(x), \quad n \in \mathbb{N}, \tag{30}$$

indicating that $\{P_n\}_0^\infty$ is a sequence of monic orthogonal polynomials. The $n$th partial numerators $\{\lambda_1\,P_{n-1}^{(1)}\}_0^\infty$ also satisfy a recurrence relation, namely,

$$P_{-1}^{(1)} \stackrel{\text{def}}{=} 0 \quad \& \quad P_0^{(1)} \stackrel{\text{def}}{=} 1 \quad \&$$

$$P_n^{(1)}(x) = (x - c_{n+1})\,P_{n-1}^{(1)}(x) - \lambda_{n+1}\,P_{n-2}^{(1)}(x), \quad n \in \mathbb{N}. \tag{31}$$

The polynomials $\{P_n^{(1)}\}_0^\infty$ are also orthogonal with respect to an appropriate measure, and are called the monic numerator polynomials or associated polynomials to $\{P_n\}_0^\infty$. Independently from the J-fraction, (31) further suggests the introduction of $\{P_n^{(k)}\}_{n=0}^\infty$ for $k \geqslant 2$, $k \in \mathbb{N}$, to be defined in the following manner

$$P_{-1}^{(k)} \stackrel{\text{def}}{=} 0 \quad \& \quad P_0^{(k)} \stackrel{\text{def}}{=} 1 \quad \&$$

$$P_n^{(k)}(x) \stackrel{\text{def}}{=} (x - c_{n+k})\,P_{n-1}^{(k)}(x) - \lambda_{n+k}\,P_{n-2}^{(k)}(x), \quad n \in \mathbb{N}. \tag{32}$$

The $\{P_n^{(k)}\}_{n=0}^\infty$ are called monic $k$th associated polynomials to $\{P_n\}_0^\infty$. Let $\psi$ be the measure with respect to which $\{P_n\}_0^\infty$ forms an orthogonal polynomial sequence.

In the case when $\operatorname{supp}\psi = [\xi_1, \eta_1]$, the J-fraction converges locally uniformly to $F$ on $\mathbb{C} - [\xi_1, \eta_1]$, where

$$F(z) = \int_{\xi_1}^{\eta_1} \frac{d\psi(x)}{z - x}, \quad z \notin [\xi_1, \eta_1]. \tag{33}$$

The measure $\psi$ can be recovered by the Stieltjes inversion formula

$$\frac{\psi(t+) + \psi(t-)}{2} = \text{const} - \frac{1}{\pi} \lim_{y \to 0+} \int_s^t \Im F(x + i\, y)\, dx\,, \tag{34}$$

(cf. [2, p. 89–90]).

Similarly to the notion on the real line, one can introduce the associated measures on the unit circle, and establish related formulas using the results of the previous sections. Let us recall that the connection between the Schur functions and the orthogonal polynomials on the unit circle was established by their common Carathéodory function ($\stackrel{\text{def}}{=} F$)

$$\frac{1 + z\, f(z)}{1 - z\, f(z)} \stackrel{\text{def}}{=} F(z) \stackrel{\text{def}}{=} \frac{1}{2\pi} \int_0^{2\pi} \frac{e^{i\theta} + z}{e^{i\theta} - z}\, d\sigma(\theta)\,. \tag{35}$$

Let us emphasize that we can start either with $\sigma$ or $f$, and determine the other one. Let

$$F_k(z) \stackrel{\text{def}}{=} \frac{1 + z\, f_k(z)}{1 - z\, f_k(z)}\,, \qquad k \in \mathbb{Z}^+\,, \tag{36}$$

where $f_k$ is determined by (3). Then, since $F_k$ is still a Carathéodory function, a measure $\sigma_k$ exists, such that

$$F_k(z) \stackrel{\text{def}}{=} \frac{1}{2\pi} \int_0^{2\pi} \frac{e^{i\theta} + z}{e^{i\theta} - z}\, d\sigma_k(\theta)\,. \tag{37}$$

We call $\sigma_k$ the $k$th associated measure to $\sigma$ ($\sigma_0 \stackrel{\text{def}}{=} \sigma$). Now one can consider the reflection coefficients $\{a_n(\sigma_k)\}_{n=0}^\infty$ of the sequence of orthogonal polynomials on the unit circle generated by $\sigma_k$. On the other hand, following [6, Example 5, p. 462], denote $\gamma_n(f_k)$ the $n$th Schur parameter of the sequence of Schur functions starting with $f_k$ (one can think of $k$ as being fixed). Then by Geronimus' theorem

$$a_n(\sigma_k) = \gamma_n(f_k) = \gamma_{n+k}(f) = a_{n+k}(\sigma)\,, \tag{38}$$

which motivates the following definition.

**Definition 11.** Let $\{a_n\}_0^\infty$ be the sequence of reflection coefficients corresponding to $\{\Phi_n\}_0^\infty$, and $k \in \mathbb{N}$. Then the polynomials $\{\Phi_n^{(k)}\}_{n=0}^\infty$ and $\{\Psi_n^{(k)}\}_{n=0}^\infty$ — as determined by the following matrix recursion — are called monic $k$th associated polynomials, and monic second kind $k$th associated polynomials, respectively

$$\Phi_0^{(k)} \stackrel{\text{def}}{=} 1 \quad \& \quad \Psi_0^{(k)} \stackrel{\text{def}}{=} 1 \quad \&$$

$$\begin{pmatrix} \Phi_{n+1}^{(k)} & \Psi_{n+1}^{(k)} \\ \Phi_{n+1}^{(k)*} & -\Psi_{n+1}^{(k)*} \end{pmatrix} \stackrel{\text{def}}{=} \begin{pmatrix} z & -\overline{a}_{n+k} \\ -z\, a_{n+k} & 1 \end{pmatrix} \begin{pmatrix} \Phi_n^{(k)} & \Psi_n^{(k)} \\ \Phi_n^{(k)*} & -\Psi_n^{(k)*} \end{pmatrix}, \qquad n \in \mathbb{Z}^+\,. \tag{39}$$

This way we are led to state the main result of the section.

**Theorem 12.** *The monic $k$th associated polynomials are orthogonal with respect to the measure $\sigma_k$ having its Carathéodory function representation*

$$F_k = \frac{F(\Phi_k + \Phi_k^*) + \Psi_k - \Psi_k^*}{F(\Phi_k - \Phi_k^*) + \Psi_k + \Psi_k^*}. \tag{40}$$

*The orthogonality measure of $\{\Psi_n^{(k)}\}_{n=0}^\infty$ has the representing Carathéodory function $1/F_k$. The moments $\{c_n(\sigma_k)\}_{n=0}^\infty$ of the measure $\sigma_k$ can be determined by*

$$c_0(\sigma_k) \stackrel{\text{def}}{=} 1 \quad \& \quad c_n(\sigma_k) \stackrel{\text{def}}{=} \frac{1}{2\pi}\int_0^{2\pi} e^{-in\theta}\,d\sigma_k(\theta) = \frac{F_k^{(n)}(0)}{2\,n!}, \quad n \in \mathbb{N}. \tag{41}$$

*Proof.* To verify (40) $F_k$ can be determined by (36), (4), (2), and (22).

For the second kind orthogonal polynomials we can always replace the initial Schur function by its negative, and the initial Carathéodory function by its reciprocal.

To show (41) we differentiate both sides of (37) $n$-times. $\square$

*Remark 13.* Let us notice that starting with the monic second kind polynomials $\{\Psi_n\}_0^\infty$, and constructing the monic $k$th associated polynomials accordingly, we are not recovering $\{\Psi_n^{(k)}\}_{n=0}^\infty$. Indeed, the Carathéodory function for the first described sequence is

$$\frac{F(\Psi_k + \Psi_k^*) + \Phi_k - \Phi_k^*}{F(\Psi_k - \Psi_k^*) + \Phi_k + \Phi_k^*},$$

which is different from $1/F_k$ unless $k = 0$.

It is a simple consequence of the definition, that the monic $k$th associated polynomials can be computed explicitly.

**Theorem 14.** *Let $k \in \mathbb{N}$, and $n \in \mathbb{Z}^+$. Then*

$$\left[2z^k \prod_{j=0}^{k-1}\left(1 - |a_j|^2\right)\right] \begin{pmatrix} \Phi_n^{(k)} & \Psi_n^{(k)} \\ \Phi_n^{(k)*} & -\Psi_n^{(k)*} \end{pmatrix} = \begin{pmatrix} \Phi_{n+k} & \Psi_{n+k} \\ \Phi_{n+k}^* & -\Psi_{n+k}^* \end{pmatrix} \begin{pmatrix} \Psi_k^* + \Psi_k & \Psi_k^* - \Psi_k \\ \Phi_k^* - \Phi_k & \Phi_k^* + \Phi_k \end{pmatrix}. \tag{42}$$

*Proof.* Using (39) and (15) iteratively, we get

$$\begin{pmatrix} \Phi_n^{(k)} & \Psi_n^{(k)} \\ \Phi_n^{(k)*} & -\Psi_n^{(k)*} \end{pmatrix} = \prod_{j=1}^n \begin{pmatrix} z & -\overline{a}_{n+k-j} \\ -z\,\overline{a}_{n+k-j} & 1 \end{pmatrix} \begin{pmatrix} 1 & 1 \\ 1 & -1 \end{pmatrix}$$

$$= \begin{pmatrix} \Phi_{n+k} & \Psi_{n+k} \\ \Phi_{n+k}^* & -\Psi_{n+k}^* \end{pmatrix} \begin{pmatrix} \Phi_k & \Psi_k \\ \Phi_k^* & -\Psi_k^* \end{pmatrix}^{-1} \begin{pmatrix} 1 & 1 \\ 1 & -1 \end{pmatrix},$$

form which (42) follows by virtue of

$$\Phi_k \Psi_k^* + \Psi_k \Phi_k^* = 2\,z^k\,\frac{1}{\kappa_k^2} = 2\,z^k \prod_{j=0}^{k-1}\left(1 - |a_j|^2\right),$$

(cf. formulas (24) and (17)). $\square$

*Remark 15.* Theorems 12 and 14 were also proved by F. Peherstorfer in [10, Theorem 3.1] (cf. [11, Theorem 2.3, p. 105]).

Department of Mathematics, The Ohio State University, 231 West Eighteenth Avenue, Columbus, Ohio 43210–1174, U.S.A.
*E-mail address*: pinter@math.ohio-state.edu

Department of Mathematics, The Ohio State University, 231 West Eighteenth Avenue, Columbus, Ohio 43210–1174, U.S.A.
*E-mail address*: nevai@math.ohio-state.edu — *WWW home*: http://www.math.ohio-state.edu/~nevai